\documentclass[12pt]{article}
\usepackage{amsmath}
\usepackage{amssymb}
\usepackage{amsfonts}
\usepackage{graphics}

\newtheorem{thm}{Theorem}[section]
\newtheorem{lemma}[thm]{Lemma}

\newtheorem{prop}[thm]{Proposition}

\newenvironment{remark}{\par\medskip\noindent{\bf Remark.\ }}{\par\smallskip}

\newcommand{\be}{\begin{equation}}
\newcommand{\ee}{\end{equation}}
\newcommand{\openbox}{\leavevmode
  \hbox to8pt{\hfil\vrule\vbox to6pt{\hrule width6pt\vfil\hrule}\vrule}}

\newcommand{\ve}[1]{\mathbf{#1}}
\newcommand{\qed}{\hbox to5pt{ } \hfill \openbox\bigskip\medskip}

\newcommand{\Fp}{\mathbb F _p}

\newcommand{\Fxv}{\F\left[\ve x\right]}

\newcommand{\cF}{\mbox{$\cal F$}}
\newcommand{\cG}{\mbox{$\cal G$}}
\newcommand{\cH}{\mbox{$\cal H$}}

\newcommand{\cM}{\mbox{$\cal M$}}

\newcommand{\N}{\mathbb N}
\newcommand{\Z}{\mathbb Z}

\newcommand{\F}{\mathbb F}

\newcommand{\rank}{\mathop\textup{rank}}
\title{$L$-balancing families}
\author{G\'abor Heged\H{u}s
\\{\normalsize  \'Obuda University}
\\{\normalsize B\'ecsi \'ut 96, Budapest, Hungary, H-1037}
\\{\normalsize hegedus.gabor@nik.uni-obuda.hu}
}

\begin{document}
\maketitle

\begin{abstract}
P. Hrube\v s, S. Natarajan Ramamoorthy, A. Rao and A. Yehudayoff proved  the following result: 

Let $p$ be a prime and let $f\in \mathbb F _p[x_1,\ldots,x_{2p}]$ be a polynomial. Suppose that $f(\mathbf{v_F})=0$ for each $F\subseteq [2p]$, where $|F|=p$  and that $f(\mathbf{0})\neq 0$. Then  $\mbox{deg}(f)\geq p$.

We prove here the following generalization of their result.

Let $p$ be a prime
and $q=p^\alpha>1$, $\alpha\geq 1$. Let $n>0$ be a positive integer and $q-1\leq d\leq n-q+1$ be an integer.
Let $\mathbb F$ be a field of characteristic $p$.
Suppose that $f(\mathbf{v_F})=0$ for each $F\subseteq [n]$, where $|F|=d$ and $\mbox{deg}(f)\leq q-1$. Then $f(\mathbf{v_F})=0$ for each $F\subseteq [n]$, where $|F|\equiv d \mbox{ (mod }q)$.

Let $t=2d$ be an even number and $L\subseteq [d-1]$ be a given subset. We say  that $\mbox{$\cal F$}\subseteq
2^{[t]}$ is an {\em $L$-balancing family} if for each $F\subseteq [t]$, where $|F|=d$ there exists a $G\subseteq [n]$ such that $|F\cap G|\in L$. 

We give a general upper bound for the size of an $L$-balancing family.
\end{abstract}
\medskip
\section{Introduction}


First we introduce some notations. 

Let $n$ be a positive integer and let $[n]$ stand for the set $\{1,2,
\ldots, n\}$. The family of all subsets of $[n]$ is denoted by $2^{[n]}$. 
For an integer $0\leq d\leq n$ we denote by
${[n] \choose d}$ the family of all  $d$ element subsets of $[n]$,
 and ${[n] \choose \leq d}={[n] \choose 0}\cup\ldots\cup{[n] \choose d}$
the subsets of size at most $d$.

Let $\F$ be a field. 
$\F[x_1, \ldots, x_n]$ denotes  the
ring of polynomials in variables $x_1, \ldots, x_n$ over $\F$.
Let $S=\F[x_1,\ldots,x_n]$. In this paper $\F$ will be a finite prime
field $\Fp$. 

In the following  $v_F\in \{0,1\}^n$ denotes the
characteristic vector of a set
$F \subseteq [n]$.
For a family of subsets $\cF \subseteq 2^{[n]}$, let
$$
V(\cF) = \{v_F : F \in \cF\} \subseteq \{0,1\}^n \subseteq \F^n.$$
It is natural to consider the ideal $I(V(\cF))$:
$$ 
I(V(\cF)):=\{f\in S:~f(v)=0 \mbox{ whenever } v\in V(\cF)\}. 
$$

Denote by $\F[x_1,\ldots,x_n]_{\leq s}$ the vector space of all polynomials
over $\F$ with degree at most $s$.

Let $I$ be an ideal of the ring $S=\F[x_1,\ldots ,x_n]$. Let $h_{S/I}(m)$ denote the dimension over $\F$ of the factor-space
 $\F[x_1,\ldots,x_n]_{\leq m}/(I\cap\F[x_1,\ldots,x_n]_{\leq m})$ (see
\cite[Section 9.3]{CLS}). The {\em Hilbert
function} of the algebra $S/I$ is the sequence $h_{S/I}(0), h_{S/I}(1),
 \ldots $.  

It is easy to verify that in the special case when $I=I(V(\cF))$ for some set system $\cF\subseteq
2^{[n]}$,
the number $h_{\cF}(m):=h_{S/I}(m)$ is the dimension of the space of functions
from
$V(\cF)$ to $\F$ which can be represented as polynomials of degree at
most $m$.

Let $p$ be a prime
and $n>1$,  $0\leq d\leq n$ be integers. Let $q=p^\alpha$, $\alpha\geq 1$.
Define the family of sets
$$
\cF(d,q)=\{K \subseteq [n]:~|K| \equiv d \mbox{ (mod }q)\}.
$$

I proved the following result in \cite{H}. 
\begin{lemma}\label{Hlemma}
Let $p$ be a prime and let $f\in \Fp[x_1,\ldots,x_{4p}]$ be a polynomial. Suppose that $f\in I(V{[4p] \choose 2p })$ and that $f\notin I(V{[4p] \choose 3p })$. Then  $deg(f)\geq p$.
\end{lemma}

My proof used a combination of Gr\"obner basis methods and linear algebra. Srinivasan gave a simpler proof which combined Fermat's little Theorem with  linear algebra (see \cite{AKV}).
Alon found a third proof based on the Combinatorial Nullstellensatz (see \cite{A}).


P. Hrube\v s, S. Natarajan Ramamoorthy, A. Rao and A. Yehudayoff proved a similar result to our Lemma \ref{Hlemma}.

\begin{lemma}\label{Hrubes}
Let $p$ be a prime and let $f\in \Fp[x_1,\ldots,x_{2p}]$ be a polynomial. Suppose that $f\in I(V{[2p] \choose p })$ and that $f(\ve 0)\neq 0$. Then  $\deg(f)\geq p$.
\end{lemma}

Let $m$ be a positive integer and and $n$ be an positive even integer. We say that a proper non-empty subsets $S_1, \ldots ,S_m$ are a {\em balancing set of family} if for every $X\in {[n]\choose n/2}$ there is an index $i\in [k]$ such that $|S_i\cap X|=|S_i|/2$.

Let $n$ be an positive even integer. 
Let $B(n)$ denote the minimum $k$ for which a balancing set family of size $k$ exists. P. Hrube\v s, S. Natarajan Ramamoorthy, A. Rao and A. Yehudayoff applied Lemma \ref{Hrubes} to give a lower bound for $B(n)$ if $n$ is an even positive integer. They used their  lower bounds to prove lower bounds on depth-2 majority and threshold circuits that compute the majority and the weighted threshold functions.  

We give here a simple  generalization of Lemma \ref{Hlemma} and Lemma \ref{Hrubes}.  Our proof method is based on the following technical result.
\begin{thm} \label{main}
Let $\F$ be an arbitrary  field.  Let $\cF\subseteq \cG\subseteq  {\F}^n$ be affine subsets. Suppose that  $m\geq 0$ is an integer such that $h_{S/I(\cF)}(m)=h_{S/I(\cG)}(m)$. Then $I(\cF)_{\leq m}=I(\cG)_{\leq m}$. 
\end{thm}

\begin{remark} 
We can easily prove  a weaker version of Combinatorial Nullstellensatz using Theorem \ref{main}. Namely let $T_i \subseteq  {\F}$ be finite subsets and let   $\cG:=\prod_{i=1}^n T_i\subseteq  {\F}^n$, where $t_i=|T_i|\geq 2$ are finite for each $i$. Let $\ve w=(w_1,\ldots ,w_n)\in \cG$ be a fixed element. Define 
$\cF:=\cG\setminus \{\ve w\}$. Let $f\in I(\cF)$ be a polynomial such that $\deg(f)\leq \sum_i t_i-n-1$. Let $m:=\sum_i t_i-n-1$. Then it is easy to check that $h_{S/I(\cF)}(m)=h_{S/I(\cG)}(m)=(\prod_i t_i)-1$. Hence  $f\in I(\cG)$. 
\end{remark} 

One of our main result is the following generalization of Lemma \ref{Hlemma} and Lemma \ref{Hrubes}.  

\begin{thm} \label{main2}
Let $p$ be a prime
and $q=p^\alpha>1$, $\alpha\geq 1$. Let $n>0$ be a positive integer and $q-1\leq d\leq n-q+1$ be an integer.
Let $\F$ be a field of characteristic $p$.
 Suppose that $f\in I(V{[n] \choose d })$ and $\mbox{deg}(f)\leq q-1$. Then $f\in I(V(\cF(d,q)))$. 
\end{thm}

We define now $L$-balancing families.
Let $n=2d$ be an even number and $L\subseteq [d-1]$ be a given subset. We say  that $\cF\subseteq
2^{[n]}$ is an {\em $L$-balancing family} if for each $F\in {[n]\choose d}$ there exists a $G\subseteq [n]$ such that $|F\cap G|\in L$. 

We prove the following general upper bound for the size of an $L$-balancing family. Our proof is based completely on Lemma \ref{Hrubes}.

\begin{thm} \label{main3}
Let $p$ be a prime. Let $n:=2p$ and $L\subseteq [p-1]$ be a given subset. Define $s:=|L|$. Let $\cF\subseteq
2^{[n]}$ be an  $L$-balancing family. Then
$$
m:=|\cF|\geq \frac{n}{2s}.
$$
\end{thm}

We prove our  results in Section 2.

\section{Proofs}

{\bf Proof of Theorem \ref{main}:}
It follows from the definition of the Hilbert function that
$$
h_{\cF}(m)=dim(S_{\leq m})-dim(I(\cF)_{\leq m})
$$
and 
$$
h_{\cG}(m)=dim(S_{\leq m})-dim(I(\cG)_{\leq m}).
$$
Since $h_{\cF}(m)=h_{\cG}(m)$, hence
$dim(I(\cF)_{\leq m})=dim(I(\cG)_{\leq m}).$

But $\cF\subseteq \cG$ implies that  $I(\cG)_{\leq m}\subseteq I(\cF)_{\leq m}$, consequently $I(\cF)_{\leq m}=I(\cG)_{\leq m}$. \qed

{\bf Proof of Theorem \ref{main2}:}
We gave an alternative proof in \cite{HR} Corollary 3.1  using Gr\"obner basis theory for Wilson's theorem about the Hilbert function of complete uniform families.
 
\begin{thm}{\em (Wilson, \cite{W}) } \label{hthm_unif}
Let $0\leq d\leq n$, $0\leq m\leq
\min\{d,n-d\}$, and $\F$ be an arbitrary field. Then we have
\begin{equation}
h_{{[n] \choose d}}(m)={n \choose m}.
\end{equation}
\end{thm}

We determined the Hilbert function of the set system  $\cF(d,q)$ in \cite{FHR} Corollary 4.5. 

\begin{thm}\label{Hilbert-thm_modp}
et $p$ be a prime
and $q=p^\alpha>1$, $\alpha\geq 1$.  Let $\F$ be a field of characteristic $p$.Let  $n>0$,  $0\leq d\leq n$ be integers
and define $r=\min\{d,n-d\}$. Let $h_{\cF(d,q)}(m)$ denote the Hilbert function 
of $\Fxv/I(V(\cF(d,q)))$. Then
\[
h_{\cF(d,q)}(m) =\sum_{i=0}^{\left\lfloor\frac mq\right\rfloor}
\binom n{m-iq}
\]
if $0\le m\le r$, and
\[
h_{\cF(d,q)}(m) =\sum_{i=-\left\lfloor\frac rq\right\rfloor}^{
\left\lfloor\frac{n-r}q\right\rfloor}
\binom n{r+iq}-
\sum_{i=1}^{\left\lfloor\frac{n-m}q\right\rfloor}
\binom n{m+iq}
\]
if $m> r$.
\end{thm}

Let $q-1\leq d\leq n-q+1$ be an integer. Suppose that $f\in I(V{[n] \choose d })_{\leq q-1}$. It follows from Theorem \ref{hthm_unif} and  Theorem \ref{Hilbert-thm_modp} that 
$h_{{[n] \choose d}}(q-1)=h_{\cF(d,q)}(q-1)={n\choose q-1}$.

Hence Theorem \ref{main} gives us that $f\in I(V(\cF(d,q)))_{\leq q-1}$. \qed

{\bf Proof of Theorem \ref{main3}:}

Let $\cF=\{F_1, \ldots ,F_m\}$ be an $L$-balancing family and let $\ve v_i:=\ve v_{F_i}$ denote the characteristic vector of $F_i$ for each $1\leq i\leq m$.

Consider the polynomial
$$
P(\ve x):=\prod_{i=1}^m \prod_{\ell \in L} (\ve x \cdot \ve v_i -\ell)\in \Fp [\ve x].
$$
Here $\cdot$ denotes the usual scalar product. Clearly $\mbox{deg}(P)\leq ms$. 

Then $P(\ve 0)=(\prod_{\ell \in L} \ell)^m\neq 0$. 
On the other hand, $P(\ve v_G)=0$ for each $G\in {[n] \choose p}$, because $\cF$ is an $L$-balancing family. 

Hence it follows from Lemma \ref{Hrubes}   that $\mbox{deg}(P)\geq p$ and we get that $p\leq ms$. \qed


\end{document}